\documentclass[11pt]{amsart}
\usepackage{bbm}
\usepackage{amssymb, amsthm, amsmath, amscd}
\usepackage{enumerate}
\usepackage[usenames,dvipsnames]{color}
\usepackage{amsfonts}
\usepackage{mathrsfs}
\usepackage{textcomp}
\usepackage[all]{xy}

\newtheorem{thm}{Theorem}[section]
\newtheorem{cor}[thm]{Corollary}
\newtheorem{lem}[thm]{Lemma}
\newtheorem{prop}[thm]{Proposition}
\theoremstyle{definition}
\newtheorem{defn}[thm]{Definition}
\theoremstyle{remark}
\newtheorem{rem}[thm]{Remark}

\numberwithin{equation}{section}

\title[Some Permanence properties for crossed products]{Some Permanence properties for crossed products by compact group actions with the tracial Rokhlin property}

\author{Haotian Tian}
\address{School of Mathematical Sciences, Tongji University, Shanghai 200092, China}
\email{2011284@tongji.edu.cn}

\author{Xiaochun Fang}
\address{School of Mathematical Sciences, Tongji University, Shanghai 200092, China}
\email{xfang@tongji.edu.cn}

\date{Month, Day, Year}


\date{\today}

\begin{document}
	\begin{abstract}
		In this paper, we give some properties of the fixed point algebra and the crossed product of a unital separable simple infinite dimensional C*-algebra by an action of a second-countable compact group with the tracial Rokhlin property with comparison that could be deduced from the properties of its original algebra: (1) stable rank one; (2) real rank zero; (3) $\beta$-comparison; (4) Winter's $n$-comparison; (5) $m$-almost divisible; (6) weakly ($m$,$n$)-divisible.
	\end{abstract}
	
	\maketitle
	
	\section{Introduction}
	The Rokhlin property for the case of a single automorphism was originally introduced for von Neumann algebras by Connes in \cite{AC1975}. Later, the Rokhlin property for finite group actions on C*-algebras first appeared in the work of Herman and Jones in \cite{HJ1982} and \cite{HJ1983}. However, the finite group acitions with the Rokhlin property are rare. Phillips, in \cite{NCP2011}, introduced the tracial Rokhlin property for finite group actions on unital simple C*-algebras. In \cite{HW2007}, Hirshberg and Winter also introduced the Rokhlin property for second-countable compact group actions on unital simple C*-algebras. Since then, crossed products by compact group actions with the Rokhlin property have been studied by several authors (see \cite{HW2007,EG2017,EG2019}). For the non-unital case, Santiago and Gardella studyed the Rokhlin property for finite group actions on non-unital simple C*-algebras in \cite{LS2015} and \cite{GS2016}. Forough and Golestani studyed the tracial Rokhlin property for finite group actions on non-unital simple C*-algebras in \cite{FG2020}. More recently, Mohammadkarimi and Phillips studied the tracial Rokhlin property with comparison for compact group actions and proved that the crossed product of a unital separable simple infinite dimensional C*-algebra with tracial rank zero by an action of a second-countable compact group with the tracial Rokhlin property with comparison has again tracial rank zero in \cite{MP2022} and some other permanence properties. Moreover, they gave some examples of compact group actions with the tracial Rokhlin property with comparison.
	
	In \cite{FF2009}, Fan and Fang proved that for a unital separable simple infinite dimensional C*-algebra $A$, a finite group $G$ and an action $\alpha \colon G\to \mathrm{Aut}(A)$ with the tracial Rokhlin property, if $A$ has stable rank one then the crossed product $A\rtimes_\alpha G$ has stable rank one, and if $A$ has real rank zero then the crossed product $A\rtimes_\alpha G$ has real rank zero. In \cite{EG2017}, Gardella proved that for a $\sigma$-unital C*-algebra $A$, a second-countable compact group $G$ and an action $\alpha \colon G\to
	\mathrm{Aut}(A)$ with the Rokhlin property, if $A$ has stable rank one then the fixed point algebra $A^\alpha$ and the crossed product $A\rtimes_\alpha G$ have stable rank one, and if $A$ has real rank zero then the fixed point algebra $A^\alpha$ and the crossed product $A\rtimes_\alpha G$ have real rank zero. In this paper, we prove that for a unital separable simple infinite dimensional C*-algebra $A$, a second-countable compact group $G$ and an action $\alpha \colon G\to \mathrm{Aut}(A)$ with the tracial tracial Rokhlin property with comparison, if $A$ has stable rank one then the fixed point algebra $A^\alpha$ and the crossed product $A\rtimes_\alpha G$ have stable rank one, and if $A$ has real rank zero then the fixed point algebra $A^\alpha$ and the crossed product $A\rtimes_\alpha G$ have real rank zero.
	
	Besides, comparison is an important property of C*-algebra. Toms and Winter conjecture that strict comparison of positive elements, finite nuclear dimension and $\mathcal{Z}$-stability are equivalent for unital separable nuclear infinite dimensional C*-algebras in \cite{WW2010}. Kirchberg and R\o rdam introduced the property of $\beta$-comparison in \cite{EM2014}. The property of $n$-comparison was introduced by Winter in \cite{WW2012}. Gardella, Hirshberg, Santiago and Vaccaro explored connections between comparison and the tracial Rokhlin property in \cite{GHS2021} and \cite{GHV2022}. In \cite{MP2022}, Mohammadkarimi and Phillips proved that the radius of comparison of the fixed point algebra of an unital separable simple infinite dimensional C*-algebra by an action of a second-countable compact group with the tracial Rokhlin property with comparison is no larger than the radius of comparison of the original algebra. In this paper, we prove that for a unital separable simple infinite dimensional C*-algebra $A$, a second-countable compact group $G$ and an action $\alpha \colon G\to \mathrm{Aut}(A)$ with the tracial tracial Rokhlin property with comparison, if $A$ has $\beta$-comparison then the fixed point algebra $A^\alpha$ and the crossed product $A\rtimes_\alpha G$ have $\beta$-comparison, and if $A$ has Winter's $n$-comparison then the fixed point algebra $A^\alpha$ and the crossed product $A\rtimes_\alpha G$ have Winter's $n$-comparison.
	
	Divisibility is also important for C*-algebra. The property of $m$-almost divisibility was introduced by Robert and Tikuisis in \cite{RT2017}. The property of weakly ($m$,$n$)-divisibility was introduced by Robert and R\o rdam in \cite{RR2013}. In this paper, we prove that for a unital separable simple infinite dimensional C*-algebra $A$, a second-countable compact group $G$ and an action $\alpha \colon G\to \mathrm{Aut}(A)$ with the tracial tracial Rokhlin property with comparison, if $A$ is $m$-almost divisible then the fixed point algebra $A^\alpha$ and the crossed product $A\rtimes_\alpha G$ are $m$-almost divisible, and if $A$ is weakly ($m$,$n$)-divisible then the fixed point algebra $A^\alpha$ and the crossed product $A\rtimes_\alpha G$ are weakly ($m$,$n$)-divisible.
	
	The regular properties considered in this paper are not at all independent. In \cite[Theorem 8.12]{APRT22} it is shown that a separable C*-algebra with stable rank one and no nonzero elementary ideal-quotients (for example, a simple, separable, stable rank one C*-algebra) has strict comparison whenever it has $\beta$-comparison for some $\beta$, or Winter's $n$-comparison for some $n$. The simple case was proved earlier in \cite{HT2020}. Following terminology of Winter \cite{WW2010}, one says that a C*-algebra is ($n$,$m$)-pure if its Cuntz semigroup has Winter's $n$-comparison and is $m$-almost divisible. By \cite[Theorem D]{APRT23}, a simple C*-algebras is pure (that is, ($0$,$0$)-pure) whenever it is ($n$,$m$)-pure for some $n$, $m$.
	
	To be precise, we get the following results.
	
	\begin{thm}
		The following properties of unital separable simple infinite dimensional C*-algebras by actions of second-countable compact groups with the tracial Rokhlin property with comparison pass from the original algebras to the fixed point algebras and the crossed products:
		\begin{enumerate}
			\item stable rank one;
			
			\item real rank zero;
			
			\item $\beta$-comparison;
			
			\item Winter's $n$-comparison;
			
			\item $m$-almost divisible;
			
			\item weakly ($m$,$n$)-divisible.
		\end{enumerate}
	\end{thm}
	
	The paper is organized as follows. Section 2 contains some preliminaries about central sequence algebras, Cuntz subequivalence and actions of second-countable compact groups with the tracial Rokhlin property with comparison. Section 3 shows the permanence of stable rank one and real rank zero. Section 4 shows the permanence of $\beta$-comparison and Winter's $n$-comparison. Section 5 shows the permanence of $m$-almost divisibility and weakly ($m$,$n$)-divisibility.
	
	\section{Preliminaries and Definitions}
	In this section, we recall some definitions and known facts about central sequence algebras, Cuntz subequivalence and the tracial Rokhlin property with comparison for second-countable compact group actions.
	
	\begin{defn}
		Let $A$ be a unital C*-algebra. We use $l^\infty(\mathbb{N},A)$ to denote the set of all bounded sequences in $A$ with the supremum norm which is a unital C*-algebra with the unit as the constant sequence $1$. Let
		\[c_0(\mathbb{N},A)=\{(a_n)_{n\in\mathbb{N}}\in l^\infty(\mathbb{N},A):\lim_{n\to\infty}\|a_n\|=0\}.\]
		It is obvious that $c_0(\mathbb{N},A)$ is a closed two-side ideal in $l^\infty(\mathbb{N},A)$, and we use the notation $A_\infty$ to denote the quotient $l^\infty(\mathbb{N},A)/c_0(\mathbb{N},A)$. Denoted by $\kappa_A\colon l^\infty(\mathbb{N},A)\to A_\infty$ the quotient map. Define $\iota\colon A\to l^\infty(\mathbb{N},A)$ by $\iota(a)=(a,a,a,\cdots)$, the constant sequence, for all $a\in A$. Identify $A$ with $\kappa_A\circ\iota(A)$. Denoted by $A_\infty\cap A'$ the relative commutant of $A$ inside of $A_\infty$.
		
		Let $\alpha\colon G\to \mathrm{Aut}(A)$ be an action of $G$ on $A$, then it induced actions of $G$ on $l^\infty(\mathbb{N},A)$ and on $A_\infty$, denoted by $\alpha^\infty$ and $\alpha_\infty$. Since for any $g\in G$,
		\[(\alpha_\infty)_g(A_\infty\cap A')\subseteq A_\infty\cap A',\]
		so we also use $\alpha_\infty$ to denote the restricted action on $A_\infty\cap A'$. These actions are not necessarily continuous when $G$ is not discrete. Therefore, we set
		\[l^\infty_\alpha(\mathbb{N},A)=\{a\in l^\infty(\mathbb{N},A):g\mapsto\alpha^\infty_g(a)\ \mathrm{is}\ \mathrm{continuous}\},\]
		and $A_{\infty,\alpha}=\kappa_A(l^\infty_\alpha(\mathbb{N},A))$. Then, $A_{\infty,\alpha}$ is invariant under $\alpha_\infty$ and the action $\alpha_\infty$ is continuous.
	\end{defn}
	
	\begin{defn}
		Let $G$ be a locally compact group, we denote the action induced by left translation of $G$ on itself by $\mathrm{Lt}\colon G\to \mathrm{Aut}(C_0(G))$.
	\end{defn}
	
	\begin{defn}
		Let $A$ be a C*-algebra, $a\in A_+$ and $\varepsilon>0$. Then we denote $f(a)$ by $(a-\varepsilon)_+$, where $f(t)=max\{0,t-\varepsilon\}$ is continuous from $[0,\infty)$ to $[0,\infty)$.
	\end{defn}
	
	The following definitions related to Cuntz comparison are from \cite{JC1978}, for more information, you can refer to \cite{GP2022} and \cite{APT2018}.
	
	\begin{defn}
		Let $A$ be a C*-algebra. Let $a,b\in (A\otimes K)_+$.
		\begin{enumerate}
			\item We say that $a$ is Cuntz subequivalent to $b$ (written $a\precsim_A b$), if there is a sequence $(r_n)_{n=1}^\infty$ in $A\otimes K$ such that $\lim\limits_{n\to\infty}\|r_n^*br_n-a\|=0$.
			
			\item We say that $a$ is Cuntz equivalent to $b$ (written $a\sim_A b$), if $a\precsim_A b$ and $b\precsim_A a$. This is an equivalence relation, we use $\langle a\rangle_A$ to denote the equivalence class of $a$. With the addition operation $\langle a\rangle_A +\langle b\rangle_A=\langle a\oplus b\rangle_A$ and the order operation $\langle a\rangle_A \leq\langle b\rangle_A$ if $a\precsim_A b$, $\mathrm{Cu}(A)=(A\otimes K)_+/\sim_A$ is an ordered semigroup which we called Cuntz semigroup. $\mathrm{W}(A)=M_\infty(A)_+/\sim_A$ is also an ordered semigroup with the same operation and order as above.
		\end{enumerate}
		If $B$ is a hereditary C*-subalgebra of $A$, and $a,b\in B_+$, then it is easy to check that $a\precsim_A b\iff a\precsim_B b$.
	\end{defn}
	
	We give the following known facts about Cuntz subequivalence. Part (1) is contained in \cite[Proposition 2.6]{KR2000} (and in a slightly different form in the earlier \cite[Proposition 2.4]{MR1992}). Part (2) is contained in \cite[Proposition 2.4]{MR1992}. Part (3) is immediate (and is \cite[Lemma 2.5(i)]{KR2000}). Part (4) is \cite[Corollary 1.6]{NCP2016}. Part (5) is \cite[Lemma 2.8(ii)]{KR2000}. Part (6) is \cite[Lemma 2.8(iii)]{KR2000}, and Part (7) is \cite[Lemma 2.9]{KR2000}.
	
	\begin{lem}\label{CSP}
		Let $A$ be a C*-algebra.
		\begin{enumerate}
			\item Let $a,b\in A_+$. Then the following are equivalent:
			\begin{enumerate}
				\item $a\precsim_A b$;
				
				\item $(a-\varepsilon)_+\precsim_A b$ for all $\varepsilon>0$;
				
				\item for every $\varepsilon>0$, there exists $\delta>0$ such that $(a-\varepsilon)_+\precsim_A(b-\delta)_+$.
			\end{enumerate}
			\item Let $a,b\in A_+$. If $a\leq b$, then we have $a\precsim_A b$.
			
			\item Let $a\in A_+$ and $\varepsilon_1,\varepsilon_2>0$. Then we have $((a-\varepsilon_1)_+-\varepsilon_2)_+=(a-(\varepsilon_1+\varepsilon_2))_+$.
			
			\item Let $\varepsilon>0$ and $\lambda\geq0$, $a,b\in A_+$. If $\|a-b\|<\varepsilon$, then we have $(a-\varepsilon-\lambda)_+\precsim_A (b-\lambda)_+$.
			
			\item Let $a,b\in A_+$, then we have $a+b\precsim_A a\oplus b$.
			
			\item Let $a,b\in A_+$ with $ab=0$, then we have $a+b\sim_A a\oplus b$.
			
			\item Let $a_1,a_2,b_1,b_2\in A_+$. If $a_1\precsim_A a_2$ and $b_1\precsim_A b_2$, then we have $a_1\oplus a_2\precsim_A b_1\oplus b_2$.
		\end{enumerate}
	\end{lem}
	
	\begin{defn}\cite[Definition 1.3]{HP2015}
		Let $G$ be a compact group, and let $A$ be a C*-algebra, $B$ be a C*-algebra. Let $\alpha\colon G\to \mathrm{Aut}(A)$ and $\gamma\colon G\to \mathrm{Aut}(B)$ be actions of $G$ on $A$ and $B$. Let $F\subseteq A$ and $S\subseteq B$ be subsets, and let $\varepsilon>0$. A completely positive contractive map $\varphi\colon A\to B$ is said to be an ($F$,$S$,$\varepsilon$)-approximately central equivariant multiplicative map if:
		\begin{enumerate}
			\item $\|\varphi(xy)-\varphi(x)\varphi(y)\|<\varepsilon$ for all $x,y\in F$.
			
			\item $\|\varphi(x)a-a\varphi(x)\|<\varepsilon$ for all $x\in F$ and all $a\in S$.
			
			\item $\sup_{g\in G}\|\varphi(\alpha_g(x))-\gamma_g(\varphi(x))\|<\varepsilon$ for all $x\in F$.
		\end{enumerate}
	\end{defn}
	
	\begin{defn}\cite[Definition 1.4]{MP2022}
		Let $A$ and $B$ be C*-algebras, and let $F\subseteq A$. A completely positive contractive map $\varphi\colon A\to B$ is said to be an ($n$,$F$,$\varepsilon$)-approximately multiplicative map if whenever $m\in\{1,2,\cdots,n\}$ and $x_1,x_2,\cdots\in F$, we have
		\[\|\varphi(x_1x_2\cdots x_m)-\varphi(x_1)\varphi(x_2)\cdots\varphi(x_m)\|<\varepsilon.\]
		If $S\subseteq B$ is also given, then $\varphi$ is said to be an ($n$,$F$,$S$,$\varepsilon$)-approximately central multiplicative map if, in addition, $\|\varphi(x)a-a\varphi(x)\|<\varepsilon$ for all $x\in F$ and all $a\in S$.
	\end{defn}
	
	Now, let us recall the notion of the tracial Rokhlin property with comparison for second-countable compact group actions defined by Mohammadkarimi and Phillips in \cite{MP2022}.
	
	\begin{defn}\cite[Definition 2.4]{MP2022}\label{DEF}
		Let $G$ be a second-countable compact group, let $A$ be a unital simple infinite dimensional C*-algebra, and let $\alpha\colon G\to \mathrm{Aut}(A)$ be an action. We say that the action $\alpha$ has the tracial Rokhlin property with comparison if for any $\varepsilon>0$,any finite set $F\subseteq A$, any finite set $S\subseteq C(G)$, any $x\in A_+$ with $\|x\|=1$, and any $y\in(A^\alpha)_+\setminus\{0\}$, there exist a projection $p\in A^\alpha$ and a unital completely positive map $\psi\colon C(G)\to pAp$ such that
		\begin{enumerate}
			\item $\psi$ is an ($F$,$S$,$\varepsilon$)-approximately central equivariant multiplicative map.
			
			\item $1-p\precsim_A x$.
			
			\item $1-p\precsim_{A^\alpha} y$.
			
			\item $1-p\precsim_{A^\alpha} p$.
			
			\item $\|pxp\|>1-\varepsilon$.
		\end{enumerate}
	\end{defn}
	
	The next theorem is the key tool for transferring properties from the original algebra to the fixed point algebra.
	
	\begin{thm}\cite[Theorem 2.17]{MP2022}\label{MTF}
		Let $G$ be a second-countable compact group, let $A$ be a unital separable simple infinite dimensional C*-algebra and let $\alpha\colon G\rightarrow \mathrm{Aut}(A)$ be an action with the tracial Rokhlin property with comparison. Then for any $\varepsilon>0$, any $n\in \mathbb{N}$, any compact subset $F_1\subseteq A$, any compact subset $F_2\subseteq A^\alpha$, any $x\in A_+$ with $\|x\|=1$, and any $y\in (A^\alpha)_+\setminus\{0\}$, there exist a projection $p\in A^\alpha$ and a unital completely positive map $\varphi\colon A\to pA^\alpha p$ such that
		\begin{enumerate}
			\item $\varphi$ is an ($n$,$F_1\cup F_2$,$\varepsilon$)-approximately multiplicative map.
			
			\item $\|pa-ap\|<\varepsilon$ for all $a\in F_1\cup F_2$.
			
			\item $\|\varphi(a)-pap\|<\varepsilon$ for all $a\in F_2$.
			
			\item $\|\varphi(a)\|\geq\|a\|-\varepsilon$ for all $a\in F_1\cup F_2$.
			
			\item $1-p\precsim_A x$.
			
			\item $1-p\precsim_{A^\alpha} y$.
			
			\item $1-p\precsim_{A^\alpha}p$.
			
			\item $\|pxp\|>1-\varepsilon$.
		\end{enumerate}
	\end{thm}
	
	The Condition $1-p\precsim_{A^\alpha} p$ is just uesd to prove
	that the algebras $A\rtimes_\alpha G$ and $A^\alpha$ are Morita equivalent (see \cite[Proposition 3.7]{MP2022}), without it, we can not transfer properties from the original algebra to the crossed product. However, the other conclusions still hold if we omit this condition. In other words, without $1-p\precsim_{A^\alpha} p$, we can still transfer properties from the original algebra to the fixed point algebra. So we consider the following version of tracial Rokhlin property for compact group actions.
	
	\begin{defn}
		Let $G$ be a second-countable compact group, let $A$ be a unital simple infinite dimensional C*-algebra, and let $\alpha\colon G\to \mathrm{Aut}(A)$ be an action. We say that the action $\alpha$ has the tracial Rokhlin property if for any $\varepsilon>0$,any finite set $F\subseteq A$, any finite set $S\subseteq C(G)$, any $x\in A_+$ with $\|x\|=1$, and any $y\in(A^\alpha)_+\setminus\{0\}$, there exist a projection $p\in A^\alpha$ and a unital completely positive map $\psi\colon C(G)\to pAp$ such that
		\begin{enumerate}
			\item $\psi$ is an ($F$,$S$,$\varepsilon$)-approximately central equivariant multiplicative map.
			
			\item $1-p\precsim_A x$.
			
			\item $1-p\precsim_{A^\alpha} y$.
			
			\item $\|pxp\|>1-\varepsilon$.
		\end{enumerate}
	\end{defn}
	
	\begin{lem}
		Let $G$ be a second-countable compact group, let $A$ be a unital separable simple infinite dimensional C*-algebra, and let $\alpha\colon G \to \mathrm{Aut}(A)$ be an action of $G$ on $A$ with the tracial Rokhlin property with comparison. Let $q$ be an $\alpha$-invariant projection in $A$. Set $B=qAq$, and denote by $\beta\colon G\to \mathrm{Aut}(B)$ the compressed action of $G$. Then $\beta$ has the tracial Rokhlin property.
		
		
		
		\begin{proof}
			Let $\varepsilon>0$, let $F\subseteq B$, and let $S\subseteq C(G)$ be finite subsets, let $x\in B_+$ with $\|x\|=1$ and let $y\in(B^\beta)_+\setminus\{0\}$. Without loss of generality, we assume that $F$ and $S$ contain only contractions. We will show that there exist a projection $p\in B^\beta$ and a unital completely positive map $\psi\colon C(G)\to pBp$ such that
			\begin{enumerate}
				\item $\psi$ is an ($F$,$S$,$\varepsilon$)-approximately central equivariant multiplicative map.
				
				\item $1-p\precsim_B x$.
				
				\item $1-p\precsim_{B^\beta} y$.
				
				\item $\|pxp\|>1-\varepsilon$.
			\end{enumerate}
			Since $\alpha$ has the tracial Rokhlin property, for $F_1=F\cup\{q\}\subseteq A$, for $S\subseteq C(G)$, for $x\in B_+\subseteq A_+$ and for $y\in(B^\beta)_+\setminus\{0\}\subseteq(A^\alpha)_+\setminus\{0\}$, there exist a projection $p_1\in A^\alpha$ and a unital completely positive map $\psi_1\colon C(G)\to p_1Ap_1$ such that
			\begin{itemize}
				\item[(5)] $\psi_1$ is an ($F_1$,$S$,$\frac{\varepsilon}{3}$)-approximately central equivariant multiplicative map.
				
				\item[(6)] $1-p_1\precsim_A x$.
				
				\item[(7)] $1-p_1\precsim_{A^\alpha} y$.
				
				\item[(8)] $\|p_1xp_1\|>1-\frac{\varepsilon}{3}$.
			\end{itemize}
			
			Set $p=qp_1q\in B^\beta$ and define a unital completely positive map $\psi\colon C(G)\to pBp$ by
			\[\psi(f)=q\psi_1(f)q\quad \mathrm{for}\ \mathrm{all}\ f\in C(G).\]
			For condition (1), for $f\in S$ and $a\in F$, we have
			\begin{align*}
				&{}\|\psi(f)a-a\psi(f)\|=\|q\psi_1(f)qa-aq\psi_1(f)q\|\\
				\leq&{}\|q\psi_1(f)qa-\psi_1(f)qa\|+\|\psi_1(f)qa-aq\psi_1(f)\|+\|aq\psi_1(f)-aq\psi_1(f)q\|\\
				<&{}\frac{\varepsilon}{3}+\frac{\varepsilon}{3}+\frac{\varepsilon}{3}=\varepsilon.
			\end{align*}
			For $f\in S$ and $g\in G$, we have
			\begin{align*}
				&{}\|\psi(\mathrm{Lt}_g(f))-\beta_g(\psi(f))\|=\|q\psi_1(\mathrm{Lt}_g(f))q-\beta_g(q\psi_1(f)q)\|\\
				=&{}\|q(\psi_1(\mathrm{Lt}_g(f))-\alpha_g(\psi_1(f)))q\|\\
				\leq&{}\|\psi_1(\mathrm{Lt}_g(f))-\alpha_g(\psi_1(f))\|<\frac{\varepsilon}{3}.
			\end{align*}
			Thus for $f\in S$, we have
			\[\sup_{g\in G}\|\psi(\mathrm{Lt}_g(f))-\beta_g(\psi(f))\|\leq\frac{\varepsilon}{3}<\varepsilon.\]
			For $f_1,f_2\in S$, we have
			\begin{align*}
				&{}\|\psi(f_1f_2)-\psi(f_1)\psi(f_2)\|=\|q\psi_1(f_1f_2)q-q\psi_1(f_1)q\psi_1(f_2)q\|\\
				\leq&{}\|q\psi_1(f_1f_2)q-q\psi_1(f_1)\psi_1(f_2)q\|+\|q\psi_1(f_1)\psi_1(f_2)q-q\psi_1(f_1)q\psi_1(f_2)q\|\\
				<&{}\frac{\varepsilon}{3}+\frac{\varepsilon}{3}=\frac{2\varepsilon}{3}<\varepsilon.
			\end{align*}
			Therefore, conditon (1) holds.
			
			For condition (2), note that $1_B-p=q(1-p_1)q\precsim_A1-p_1\precsim_A x$ and $B$ is a hereditary $C^*$-algebra of $A$, it follows that $1_B-p\precsim_B x$.
			
			For condition (3), note that $1_B-p=q(1-p_1)q\precsim_{A^\alpha}1-p_1\precsim_{A^\alpha} y$ and $B^\beta$ is a hereditary $C^*$-algebra of $A^\alpha$, it follows that $1_B-p\precsim_{B^\beta} y$.
			
			For condition (4), we have
			\[\|pxp-x\|=\|qp_1qxqp_1q-x\|=\|qp_1xp_1q-qxq\|\leq\|p_1xp_1-x\|<\frac{\varepsilon}{3}.\]
			Thus, $\|pxp\|>1-\frac{\varepsilon}{3}>1-\varepsilon$.
		\end{proof}
	\end{lem}
	
	\begin{rem}
		In general, the tracial Rokhlin property with comparison usually can not pass to corners since we can not always get $q(1-p)q\precsim_{B^\beta} qpq$ from $1-p\precsim_{A^\alpha} p$ for any $\alpha$-invariant projection $q$ in $A$ and $B=qAq$.
	\end{rem}

	\begin{thm}($\mathrm{cf}$.\cite[Theorem 2.17]{MP2022})\label{HMTF}
		Let $G$ be a second-countable compact group, let $A$ be a unital separable simple infinite dimensional C*-algebra and let $\alpha\colon G\to \mathrm{Aut}(A)$ be an action with the tracial Rokhlin property. Then for any $\varepsilon>0$, any $n\in \mathbb{N}$, any compact subset $F_1\subseteq A$, any compact subset $F_2\subseteq A^\alpha$, any $x\in A_+$ with $\|x\|=1$, and any $y\in (A^\alpha)_+\setminus\{0\}$, there exist a projection $p\in A^\alpha$ and a unital completely positive map $\varphi\colon A\to pA^\alpha p$ such that
		\begin{enumerate}
			\item $\varphi$ is an ($n$,$F_1\cup F_2$,$\varepsilon$)-approximately multiplicative map.
			
			\item $\|pa-ap\|<\varepsilon$ for all $a\in F_1\cup F_2$.
			
			\item $\|\varphi(a)-pap\|<\varepsilon$ for all $a\in F_2$.
			
			\item $\|\varphi(a)\|\geq\|a\|-\varepsilon$ for all $a\in F_1\cup F_2$.
			
			\item $1-p\precsim_A x$.
			
			\item $1-p\precsim_{A^\alpha} y$.
			
			\item $\|pxp\|>1-\varepsilon$.
		\end{enumerate}
		\begin{proof}
			The proof is essentially the same as that of \cite[Theorem 2.17]{MP2022}, we only omit the condition $1-p\precsim_{A^\alpha}p$.
		\end{proof}
	\end{thm}
	
	\begin{thm}\cite[Theorem 3.9, Corollary 3.10]{MP2022}\label{CPS}
		Let $G$ be a second-countable compact group, let $A$ be a unital separable simple infinite dimensional C*-algebra and let $\alpha\colon G\to \mathrm{Aut}(A)$ be an action with the tracial Rokhlin property with comparison. Then the crossed product $A\rtimes_\alpha G$ is simple. Moreover, the algebras $A\rtimes_\alpha G$ and $A^\alpha$ are Morita equivalent and stably isomorphic.
	\end{thm}
	
	\section{Stable rank one and Real rank zero}
	In this section, we prove that the fixed point algebra and the crossed product of a unital separable simple infinite dimensional C*-algebra with stable rank one (respectively, real rank zero) by an action of a second-countable compact group with the tracial Rokhlin property with comparison have again stable rank one (respectively, real rank zero). First, we introduce the definitions of stable rank one and real rank zero.
	\begin{defn}
		Let $A$ be a C*-algebra.
		\begin{enumerate}
			\item If $A$ is unital, it is said to have stable rank one if the set of invertible elements is dense in $A$, written $\mathrm{tsr}(A)=1$. If $A$ is not unital, it is said to have stable rank one if for its unitization $\tilde{A}$, we have $\mathrm{tsr}(\tilde{A})=1$.
			
			\item If $A$ is unital, it is said to have real rank zero if the set of invertible self-adjoint elements is dense in $A_{sa}$, written $\mathrm{RR}(A)=0$. If $A$ is not unital, it is said to have real rank zero if for its unitization $\tilde{A}$, we have $\mathrm{RR}(\tilde{A})=0$.
		\end{enumerate}
	\end{defn}
	
	\begin{defn}
		A C*-algebra $A$ is said to have the Property (SP), if every nonzero hereditary C*-subalgebra of $A$ contains a nonzero projection.
	\end{defn}
	
	\begin{defn}
		Let $\varepsilon>0$. Define
		\begin{equation*}
			f_\varepsilon(t)=
			\begin{cases}
				0,&\mbox{$t\in[0,\varepsilon]$,}\\
				$linear$,&\mbox{$t\in[\varepsilon,2\varepsilon]$.}\\
				1,&\mbox{$t\in[2\varepsilon,+\infty)$.}
			\end{cases}
		\end{equation*}
	\end{defn}
	
	Now, we will show the permanence of stable rank one. Given an element $x\in A^\alpha$, we will find an invertible element in $A^\alpha$ to approximate it. First, we find a zero-divisor $c\in A^\alpha$ to approximate $x$. Next, we show that there exists a projection $p_1\in A^\alpha$ and find an invertible element $b\in A_1=(1-p_1)A(1-p_1)$ to approximate $c$. Then, we can use Theorem \ref{HMTF} to get a projection $p\in (A_1)^\alpha$ and an approximately multiplicative map $\varphi\colon A_1\to p(A_1)^\alpha p$ and the sum of two elements to approximate $c$. One of them is $\varphi(b)$ which is invertible in $p(A_1)^\alpha p$, the other is Cuntz subequivalent to $p_1$. Last, we find an invertible element in $A^\alpha$ to approximate $c$ by matrix decomposition and some perturbations.
	\begin{thm}\label{SR1}
		Let $A$ be a unital separable simple infinite dimensional C*-algebra with stable rank one. Let $\alpha\colon G\to \mathrm{Aut}(A)$ be an action of a second-countable compact group which has the tracial Rokhlin property with comparison. Then the fixed point algebra $A^\alpha$ has stable rank one.
		\begin{proof}
			If $A$ does not have the Property (SP), by \cite[Lemma 2.8]{MP2022}, $\alpha$ has the Rokhlin property. It follows that $A^\alpha$ has stable rank one (see \cite[Proposition 4.13(2)]{EG2017}).
			
			Thus we may assume that $A$ has the Property (SP). Let $a\in A^\alpha$ and let $\varepsilon\in(0,1)$. Without loss of generality, we may assume that $a$ is not invertible and $\|a\|=1$. $A$ is stably finite since it has stable rank one (see \cite[Lemma 3.5]{MAR1983}). Thus $A^\alpha$ is stably finite. Therefore $a$ is not one-sided invertible. By \cite[Lemma 3.6.9]{HL2001}, we may assume that there exists a zero-divisor $c\in A^\alpha$ such that
			\[\|a-c\|<\frac{\varepsilon}{4}.\]
			Since $A$ has the Property (SP), by \cite[Corollary 4.4]{MP2022}, so does $A^\alpha$. Then, there exists a nonzero projection $e\in A^\alpha$ such that $ce=ec=0$. Since $A^\alpha$ is not of type \uppercase\expandafter{\romannumeral1}, simple and unital (see \cite[Theorem 3.2]{MP2022} and \cite[Proposition 3.3]{MP2022}), by \cite[Lemma 1.10]{NCP2011}, there are nonzero projections $p_1,p_2\in$ Her$(e)$ such that $p_1p_2=0$ and $p_1\sim_{A^\alpha}p_2$. Set $A_1=(1-p_1)A(1-p_1)$. Since $A$ has stable rank one, it follows that $A_1$ has stable rank one. Note that $c\in A_1$, there exists an invertible element $b\in A_1$ such that 
			\[\|c-b\|<\frac{\varepsilon}{8}.\]
			Then, by Theorem \ref{HMTF}, for $F_1=\{b\}$, $F_2=\{c\}$, there exist a projection $p\in (A_1)^\alpha$ and a unital completely positive map $\varphi\colon A_1\to p(A_1)^\alpha p$ such that
			\begin{enumerate}
				\item $\varphi$ is an ($2$, $F_1\cup F_2$, $\frac{\varepsilon}{8}$)-approximately multiplicative map.
				
				\item $\|px-xp\|<\frac{\varepsilon}{8}$ for all $x\in F_1\cup F_2$.
				
				\item $\|\varphi(x)-pxp\|<\frac{\varepsilon}{8}$ for all $x\in F_2$.
				
				\item $1_{A_1}-p\precsim_{(A_1)^\alpha} p_2$.
			\end{enumerate}
			Set $c_1=pcp$ and $c_2=(1_{A_1}-p)c(1_{A_1}-p)$. Then
			\[\|c-(c_1+c_2)\|<\frac{\varepsilon}{4}.\]
			Note that
			\[\|\varphi(b)\varphi(b^{-1})-p\|=\|\varphi(b^{-1})\varphi(b)-p\|<\frac{\varepsilon}{8}<1,\]
			\[\|\varphi(c)-pcp\|<\frac{\varepsilon}{8}.\]
			Then $\varphi(b)\varphi(b^{-1})$ and $\varphi(b^{-1})\varphi(b)$ are invertible in $p(A_1)^\alpha p=p(1-p_1)A^\alpha(1-p_1)p$. It follows that $\varphi(b)$ is invertible. Thus we have
			\begin{align*}
				\|c_1-\varphi(b)\|\leq&{}\|c_1-\varphi(c)\|+\|\varphi(c)-\varphi(b)\|\\
				<&{}\frac{\varepsilon}{8}+\frac{\varepsilon}{8}=\frac{\varepsilon}{4}.
			\end{align*}
			By (4), there exists $u\in A^\alpha$ such that 
			\[u^*u=1_{A_1}-p=(1-p)(1-p_1)\ \mathrm{and}\ uu^*\leq p_1.\]
			Set $d=c_2+\frac{\varepsilon}{16}u+\frac{\varepsilon}{16}u^*+\frac{\varepsilon}{4}(p_1-uu^*)$. Therefore $d$, in a matrix form respect to $(1-p)(1-p_1)$, $p(1-p_1)$, $p_1$;
			\begin{equation*}
				\begin{pmatrix}
					c_2 &\frac{\varepsilon}{16} & 0\\
					\frac{\varepsilon}{16} & 0& 0\\
					0 & 0 & \frac{\varepsilon}{4} 
				\end{pmatrix},
			\end{equation*}
			is an invertible element in $(1-p)A^\alpha(1-p)=(1-p)(1-p_1)A^\alpha(1-p_1)(1-p)+p_1A^\alpha p_1$. Moreover, 
			\[\|d-c_2\|<\frac{\varepsilon}{4}.\]
			Hence, $\varphi(b)+d$ is invertible in $A^\alpha$. Finally,
			\begin{align*}
				\|a-(\varphi(b)+d)\|\le&{}\|a-c\|+\|c-(\varphi(b)+d)\|\\
				<&{}\frac{\varepsilon}{4}+\|c-(c_1+c_2)\|+\|(c_1+c_2)-(\varphi(b)+d)\|\\
				<&{}\frac{\varepsilon}{4}+\frac{\varepsilon}{4}+\|c_1-\varphi(b)\|+\|c_2-d\|\\
				<&{}\frac{\varepsilon}{4}+\frac{\varepsilon}{4}+\frac{\varepsilon}{4}+\frac{\varepsilon}{4}\\
				=&{}\varepsilon.
			\end{align*}
		\end{proof}
	\end{thm}
	
	\begin{cor}
		Let $A$ be a unital separable simple infinite dimensional C*-algebra with stable rank one. Let $\alpha\colon G\to \mathrm{Aut}(A)$ be an action of a second-countable compact group which has the tracial Rokhlin property with comparison. Then the crossed product $A\rtimes_\alpha G$ has stable rank one.
		\begin{proof}
			The algebra $A^\alpha$ has stable rank one by Theorem \ref{SR1}. Theorem \ref{CPS} implies that $A\rtimes_\alpha G$ is Morita equivalent to $A^\alpha$, so \cite[Theorem 3.6]{MAR1983} implies that $A\rtimes_\alpha G$ has stable rank one.
		\end{proof}
	\end{cor}
	
	Now, we will show the permanence of real rank zero. The idea of the proof is essentially the same as that of Theorem \ref{SR1}, the only difference is replacing invertible elements with invertible self-adjoint elements.
	\begin{thm}\label{RR0}
		Let $A$ be a unital separable simple infinite dimensional C*-algebra with real rank zero. Let $\alpha\colon G\to \mathrm{Aut}(A)$ be an action of a second-countable compact group which has the tracial Rokhlin property with comparison. Then the fixed point algebra $A^\alpha$ has real rank zero.
		\begin{proof}
			Since $A$ has real rank zero, we know that $A$ has the Property (SP). Let $a\in (A^\alpha)_{sa}$ and let $\varepsilon\in(0,1)$.  Without loss of generality, we may assume that $a$ is not invertible and $\|a\|=1$.
			Define a continuous function $f\in C([-1,1])$ such that $0\leq f\leq 1$ and $f(t)=0$ if $|t|\geq \frac{\varepsilon}{4}$ and $f(t)=1$ if $|t|<\frac{\varepsilon}{8}$. Since $0\in \mathrm{sp}(a)$, $f(a)\neq 0$. By \cite[Corollary 4.4]{MP2022}, there exists a nonzero projection $e\in $Her$(f(a))\subseteq A^\alpha$. Since $A^\alpha$ is not of type \uppercase\expandafter{\romannumeral1}, simple and unital (see \cite[Theorem 3.2]{MP2022} and \cite[Proposition 3.3]{MP2022}), by \cite[Lemma 1.10]{NCP2011}, there are nonzero projections $p_1,p_2\in$ Her$(e)$ such that $p_1p_2=0$ and $p_1\sim_{A^\alpha} p_2$. Let $c=f_{\frac{\varepsilon}{4}}(a)a$. Then $p_1c=cp_1=0$ and
			\[\|a-c\|<\frac{\varepsilon}{4}.\]
			Set $A_1=(1-p_1)A(1-p_1)$. Since $A$ has real rank zero, it follows that $A_1$ has real rank zero. Note that $c\in A_1$, there exists an invertible self-adjoint element $b\in A_1$ such that 
			\[\|c-b\|<\frac{\varepsilon}{8}.\]
			As in the same step in the arguement in Theorem \ref{SR1}, we get a projection $p\in (A_1)^\alpha$ and a unital completely positive map $\varphi\colon A_1\to p(A_1)^\alpha p$ such that $\varphi(b)\in(p(A_1)^\alpha p)_{sa}=(p(1-p_1)A^\alpha(1-p_1)p)_{sa}$ is invertible and $d\in ((1-p)(A^\alpha)(1-p))_{sa}=((1-p)(1-p_1)A^\alpha(1-p_1)(1-p)+p_1A^\alpha p_1)_{sa}$ which is invertible such that
			\[\|c-(\varphi(b)+d)\|<\frac{3\varepsilon}{4}.\]
			Moreover, we have $\varphi(b)+d$ is invertible in $(A^\alpha)_{sa}$ and
			\[\|a-\varphi(b)+d\|<\varepsilon.\]
		\end{proof}
	\end{thm}
	
	\begin{cor}\label{RR01}
		Let $A$ be a unital separable simple infinite dimensional C*-algebra with real rank zero. Let $\alpha\colon G\to \mathrm{Aut}(A)$ be an action of a second-countable compact group which has the tracial Rokhlin property with comparison. Then the crossed product $A\rtimes_\alpha G$ has real rank zero.
		\begin{proof}
			The algebra $A^\alpha$ has real rank zero by Theorem \ref{RR0}. Theorem \ref{CPS} implies that $A\rtimes_\alpha G$ is Morita equivalent to $A^\alpha$, so \cite[Theorem 3.8]{BP1991} implies that $A\rtimes_\alpha G$ has real rank zero.
		\end{proof}
	\end{cor}
	
	\section{$\beta$-comparison and Winter's $n$-comparison}
	In this section, we prove that the fixed point algebra and the crossed product of a unital separable simple infinite dimensional stably finite C*-algebra with $\beta$-comparison (respectively, Winter's $n$-comparison) by an action of a second-countable compact group with the tracial Rokhlin property with comparison have again $\beta$-comparison (respectively, Winter's $n$-comparison). First, we introduce the definitions of $\beta$-comparison and Winter's $n$-comparison.
	
	\begin{defn}\cite[Definition 2.2]{EM2014}
		Let $A$ be a C*-algebra and let $1\leq\beta<\infty$. We say that $A$ has $\beta$-comparison if for all $\langle x\rangle,\langle y\rangle\in \mathrm{Cu}(A)$ and all integers $k,l\geq1$ with $k>\beta l$, the inequality $k\langle x\rangle\leq l\langle y\rangle$ implies $\langle x\rangle\leq \langle y\rangle$.
	\end{defn}
	
	\begin{defn}\cite[Definition 3.1]{WW2012}
		Let $A$ be a C*-algebra and write $x<_sy$ to mean $(k+1)\langle x\rangle\leq k\langle y\rangle$ for some $k\in\mathbb{N}$. We say tha $A$ has Winter's $n$-comparison, if, whenever $\langle x\rangle,\langle y_0\rangle,\langle y_1\rangle,\cdots,\langle y_n\rangle\in \mathrm{Cu}(A)$ such that $\langle x\rangle<_s \langle y_j\rangle$ for all $j=0,1,\cdots,n$, then $\langle x\rangle\leq \langle y_0\rangle+\langle y_1\rangle+\cdots+\langle y_n\rangle$.
	\end{defn}

	The proof of the following lemma is contained in the \cite[Proposition 4.22]{MP2022} which is the same as that of \cite[Lemma 3.10]{ANP2021}.
	
	\begin{lem}\label{EP}
		Let $G$ be a second-countable compact group, let $A$ be a unital separable simple infinite dimensional C*-algebra, and let $\alpha\colon G\to \mathrm{Aut}(A)$ be an action with the tracial Rokhlin property with comparison. Then, for every $x\in(A^\alpha)_+\setminus\{0\}$, there exists $c\in(A^\alpha)_+$ such that $c\precsim_{A^\alpha} x$ and $\mathrm{sp}(c)=[0,1]$.
	\end{lem}
	
	Now, we will show the permanence of $\beta$-comparison. The key of the proof is using \cite[Proposition 4.21]{MP2022} to transfer the Cuntz subequivalence from the original algebra to the fixed point algebra. We use \cite[Lemma 2.7]{NCP2016} to work with details.
	\begin{thm}\label{BC}
		Let $A$ be a unital separable simple infinite dimensional stably finite C*-algebra with $\beta$-comparison. Let $\alpha\colon G\to \mathrm{Aut}(A)$ be an action of a second-countable compact group which has the tracial Rokhlin property with comparison. Then the fixed point algebra $A^\alpha$ has $\beta$-comparison
		\begin{proof}
			Let $a,b\in (A^\alpha\otimes K)_+$ and let $k,l\ge1$ be integers with $k>\beta l$ such that
			\[k\langle a\rangle_{A^\alpha}\leq l\langle b\rangle_{A^\alpha}.\]
			Since $\langle a\rangle_{A^\alpha}=\sup_{\varepsilon>0} \langle (a-\varepsilon)_+\rangle_{A^\alpha}$, we need to show that $(a-\varepsilon)_+\precsim_{A^\alpha}b$ for any $\varepsilon>0$.
			
			Fix $0<\varepsilon<\frac{1}{2}$. Choose $m\in\mathbb{N}$ such that
			\[\frac{m}{m+1}\frac{k}{l}>\beta.\]
			Then we have
			\[mk\langle a\rangle_{A^\alpha}\leq ml\langle b\rangle_{A^\alpha}.\]
			Therefore, by Lemma \ref{CSP}(1), there is $\delta>0$ such that
			\[mk\langle (a-\frac{\varepsilon}{3})_+\rangle_{A^\alpha}\leq ml\langle (b-\delta)_+\rangle_{A^\alpha}.\]
			
			Choose $n\in\mathbb{N}$ and $a_0,b_0\in M_n(A^\alpha)_+$ such that
			\[\|a_0-(a-\frac{\varepsilon}{3})_+\|<\frac{\varepsilon}{3}\quad \mathrm{and}\quad \|b_0-(b-\frac{\delta}{3})_+\|<\frac{\delta}{6}.\]
			Then we have
			\[(a-\varepsilon)_+\precsim_{A^\alpha}(a_0-\frac{\varepsilon}{3})_+\precsim_{A^\alpha}(a-\frac{\varepsilon}{3})_+\]
			and
			\[(b-\delta)_+\precsim_{A^\alpha}(b_0-\frac{\delta}{3})_+\precsim_{A^\alpha}(b_0-\frac{\delta}{6})_+\precsim_{A^\alpha}(b-\frac{\delta}{3})_+.\]
			Thus we have
			\[mk\langle (a_0-\frac{\varepsilon}{3})_+\rangle_{A^\alpha}\leq ml\langle (b_0-\frac{\delta}{3})_+\rangle_{A^\alpha}.\]
			
			By \cite[Proposition 4.20]{MP2022}, the action $\alpha\otimes \mathrm{id}_{M_n}\colon G\to \mathrm{Aut}(A\otimes M_n)$ also has the tracial Rokhlin property with comparison. Thus $M_n(A^\alpha)$ is not of type \uppercase\expandafter{\romannumeral1}, simple and infinite
			dimensional, by \cite[Lemma 2.7]{NCP2016}, there are positive elements $c\in M_n(A^\alpha)_+$ and $d\in M_n(A^\alpha)_+\setminus\{0\}$ such that
			\[m\langle (b_0-\frac{\delta}{3})_+\rangle_{A^\alpha}\leq(m+1)\langle c\rangle_{A^\alpha}\quad \mathrm{and} \quad \langle c\rangle_{A^\alpha}+\langle d\rangle_{A^\alpha}\leq\langle (b_0-\frac{\delta}{6})_+\rangle_{A^\alpha}.\]
			By Lemma \ref{EP}, there is $d_0\in M_n(A^\alpha)_+$ such that $d_0\precsim_{A^\alpha}d$ and $\mathrm{sp}(d_0)=[0,1]$. Thus we have 
			\[\langle c\rangle_{A^\alpha}+\langle d_0\rangle_{A^\alpha}\leq\langle (b_0-\frac{\delta}{6})_+\rangle_{A^\alpha},\]
			and
			\[mk\langle (a_0-\frac{\varepsilon}{3})_+\rangle_{A^\alpha}\leq (m+1)l\langle c \rangle_{A^\alpha}.\]
			
			Since $A^\alpha\subset A$ and $A$ has $\beta$-comparison and
			\[mk>\beta(m+1)l,\]
			we have $(a_0-\frac{\varepsilon}{3})_+\precsim_A c$. Therefore,
			\[(a-\varepsilon)_+\precsim_{A^\alpha}(a_0-\frac{\varepsilon}{3})_+\precsim_A c\oplus d_0.\]
			Note that $0$ is a limit point of $\mathrm{sp}(c+d_0)$, by \cite[Proposition 4.21]{MP2022},
			\[(a-\varepsilon)_+\precsim_{A^\alpha}(a_0-\frac{\varepsilon}{3})_+\precsim_{A^\alpha} c\oplus d_0\precsim_{A^\alpha}(b_0-\frac{\delta}{6})_+\precsim_{A^\alpha}(b-\frac{\delta}{3})_+\precsim_{A^\alpha}b.\]
		\end{proof}
	\end{thm}
	
	\begin{cor}
		Let $A$ be a unital separable simple infinite dimensional stably finite C*-algebra with $\beta$-comparison. Let $\alpha\colon G\to \mathrm{Aut}(A)$ be an action of a second-countable compact group which has the tracial Rokhlin property with comparison. Then the crossed product $A\rtimes_\alpha G$ has $\beta$-comparison
		\begin{proof}
			The algebra $A^\alpha$ has $\beta$-comparison by Theorem \ref{BC}. Theorem \ref{CPS} implies that $A\rtimes_\alpha G$ is stably isomorphic to $A^\alpha$. Since two stably isomorphic $C^*$-algebras have canonically isomorphic Cuntz semigroups, we know that $A\rtimes_\alpha G$ has $\beta$-comparison.
		\end{proof}
	\end{cor}
	
	Now, we will show the permanence of Winter's $n$-comparison. The idea of the proof is essentially the same as that of Theorem \ref{BC}. 
	\begin{thm}\label{NC}
		Let $A$ be a unital separable simple infinite dimensional stably finite C*-algebra with Winter's $n$-comparison. Let $\alpha\colon G\to\mathrm{Aut}(A)$ be an action of a second-countable compact group which has the tracial Rokhlin property with comparison. Then fixed point algebra $A^\alpha$ has Winter's $n$-comparison. 
		\begin{proof}
			Let $a,b_0,b_1,\cdots,b_n\in (A^\alpha\otimes K)_+$ and let $k_i\in\mathbb{N}$ such that
			\[(k_i+1)\langle a\rangle_{A^\alpha}\leq k_i\langle b_i\rangle_{A^\alpha}\]
			in $\mathrm{Cu}(A)$ for all $0\leq i\leq n$. Since $\langle a\rangle_{A^\alpha}=\sup_{\varepsilon>0} \langle (a-\varepsilon)_+\rangle_{A^\alpha}$, we need to show that $(a-\varepsilon)_+\precsim_{A^\alpha} b_0\oplus b_1\oplus\cdots\oplus b_n$ for any $\varepsilon>0$.
			
			Fix $0<\varepsilon<\frac{1}{2}$. Note that $k_i$ can be choosen to be the same for all $b_i$, as, with $k=(k_0+1)(k_1+1)\cdots(k_n+1)-1$, one has
			\[(k+1)\langle a\rangle_{A^\alpha}\leq k\langle b_i\rangle_{A^\alpha}\] for all $0\leq i\leq n$. Choose $m\in\mathbb{N}$ such that
			\[m-k\geq1.\]
			Then we have
			\[m(k+1)\langle a\rangle_{A^\alpha}\leq mk\langle b_i\rangle_{A^\alpha}\]
			for all $0\leq i\leq n$. Therefore, by Lemma \ref{CSP}(1), there is $\delta>0$ such that
			\[m(k+1)\langle (a-\frac{\varepsilon}{3})_+\rangle_{A^\alpha}\leq mk\langle (b_i-\delta)_+\rangle_{A^\alpha}\]
			for all $0\leq i\leq n$. 
			
			Choose $l\in\mathbb{N}$ and $a',b_i'\in M_l(A^\alpha)_+$ for $0\leq i\leq n$ such that
			\[\|a'-(a-\frac{\varepsilon}{3})_+\|<\frac{\varepsilon}{3}\quad \mathrm{and}\quad \|b_i'-(b_i-\frac{\delta}{3})_+\|<\frac{\delta}{6}.\]
			Then we have
			\[(a-\varepsilon)_+\precsim_{A^\alpha}(a'-\frac{\varepsilon}{3})_+\precsim_{A^\alpha}(a-\frac{\varepsilon}{3})_+\]
			and
			\[(b_i-\delta)_+\precsim_{A^\alpha}(b_i'-\frac{\delta}{3})_+\precsim_{A^\alpha}(b_i'-\frac{\delta}{6})_+\precsim_{A^\alpha}(b_i-\frac{\delta}{3})_+.\]
			for all $0\leq i\leq n$. Thus we have
			\[m(k+1)\langle (a'-\frac{\varepsilon}{3})_+\rangle_{A^\alpha}\leq mk\langle (b_i'-\frac{\delta}{3})_+\rangle_{A^\alpha}\]
			for all $0\leq i\leq n$. 
			
			By \cite[Proposition 4.20]{MP2022}, the action $\alpha\otimes \mathrm{id}_{M_l}\colon G\to \mathrm{Aut}(A\otimes M_l)$ also has the tracial Rokhlin property with comparison. Thus $M_l(A^\alpha)$ is not of type \uppercase\expandafter{\romannumeral1}, simple and infinite
			dimensional, by \cite[Lemma 2.7]{NCP2016}, there are positive elements $c_i\in M_l(A^\alpha)_+$ and $d_i\in M_l(A^\alpha)_+\setminus\{0\}$ such that
			\[m\langle(b_i'-\frac{\delta}{3})_+\rangle_{A^\alpha}\leq(m+1)\langle c_i\rangle_{A^\alpha}\quad \mathrm{and} \quad \langle c_i\rangle_{A^\alpha}+\langle d_i\rangle_{A^\alpha}\leq\langle (b_i'-\frac{\delta}{6})_+\rangle_{A^\alpha}\]
			for all $0\leq i\leq n$. By Lemma \ref{EP}, there are $d_i'\in M_l(A^\alpha)_+$ such that $d_i'\precsim_{A^\alpha}d_i$ and $\mathrm{sp}(d_i')=[0,1]$ for all $0\leq i\leq n$. Thus we have
			\[\langle c_i\rangle_{A^\alpha}+\langle d_i'\rangle_{A^\alpha}\leq\langle (b_i'-\frac{\delta}{6})_+\rangle_{A^\alpha},\]
			and
			\[m(k+1)\langle (a'-\frac{\varepsilon}{3})_+\rangle_{A^\alpha}\leq (m+1)k\langle c_i \rangle_{A^\alpha}\]
			for all $0\leq i\leq n$.
			
			Since $A^\alpha\subset A$ and $A$ has Winter's $n$-comparison and
			\[m-k\geq1,\]
			we have 
			\[((m+1)k+1)\langle (a'-\frac{\varepsilon}{3})_+\rangle_{A^\alpha}\leq m(k+1)\langle (a'-\frac{\varepsilon}{3})_+\rangle_{A^\alpha}\leq (m+1)k\langle c_i \rangle_{A^\alpha}.\]
			Thus,
			\[(a'-\frac{\varepsilon}{3})_+\precsim_A c_0\oplus c_1\oplus\cdots\oplus c_n.\]
			Therefore,
			\[(a-\varepsilon)_+\precsim_{A^\alpha}(a'-\frac{\varepsilon}{3})_+\precsim_A c_0\oplus d_0'\oplus c_1\oplus d_1'\oplus\cdots \oplus c_n\oplus d_n'.\]
			Note that $0$ is a limit point of $\mathrm{sp}(c_0\oplus d_0'\oplus c_1\oplus d_1'\oplus\cdots \oplus c_n\oplus d_n')$, by \cite[Proposition 4.21]{MP2022},
			\begin{align*}
				(a-\varepsilon)_+&{}\precsim_{A^\alpha} (a'-\frac{\varepsilon}{3})_+\\
				&{}\precsim_{A^\alpha}c_0\oplus d_0'\oplus c_1\oplus d_1'\oplus\cdots \oplus c_n\oplus d_n'\\
				&{}\precsim_{A^\alpha} (b_0'-\frac{\delta}{6})_+\oplus\cdots\oplus(b_n'-\frac{\delta}{6})_+\\
				&{}\precsim_{A^\alpha}(b_0-\frac{\delta}{3})_+\oplus\cdots\oplus(b_n-\frac{\delta}{3})_+\\
				&{}\precsim_{A^\alpha}b_0\oplus\cdots\oplus b_n.
			\end{align*}
		\end{proof}
	\end{thm}
	
	\begin{cor}
		Let $A$ be a unital separable simple infinite dimensional stably finite C*-algebra with Winter's $n$-comparison. Let $\alpha\colon G\to \mathrm{Aut}(A)$ be an action of a second-countable compact group which has the tracial Rokhlin property with comparison. Then the crossed product $A\rtimes_\alpha G$ has Winter's $n$-comparison
		\begin{proof}
			The algebra $A^\alpha$ has Winter's $n$-comparison by Theorem \ref{NC}. Theorem \ref{CPS} implies that $A\rtimes_\alpha G$ is stably isomorphic to $A^\alpha$. Since two stably isomorphic $C^*$-algebras have canonically isomorphic Cuntz semigroups, we know that $A\rtimes_\alpha G$ has Winter's $n$-comparison.
		\end{proof}
	\end{cor}
	
	\begin{cor}
		Let $A$ be a unital separable simple infinite dimensional stably finite C*-algebra with strict comparison. Let $\alpha\colon G\to \mathrm{Aut}(A)$ be an action of a second-countable compact group which has the tracial Rokhlin property with comparison. Then the fixed point algrebra $A^\alpha$ and the crossed product $A\rtimes_\alpha G$ have strict comparison.
		\begin{proof}
			Apply Theorem \ref{BC} when $\beta=1$ or Theorem \ref{NC} when $n=0$.
		\end{proof}
	\end{cor}
	
	\section{$m$-almost divisibility and weakly ($m$,$n$)-divisibility}
	In this section, we prove that the fixed point algebra and the crossed product of a unital separable simple infinite dimensional stably finite C*-algebra which is $m$-almost divisible (respectively, weakly ($m$,$n$)-divisible) by an action of a second-countable compact group with the tracial Rokhlin property with comparison which are again $m$-almost divisible (respectively, weakly ($m$,$n$)-divisible). First, we introduce the definitions of $m$-almost divisibility and weakly ($m$,$n$)-divisibility.
	\begin{defn}\cite[2.3]{RT2017}
		Let $m\in\mathbb{N}$. We say that the C*-algebra $A$ is $m$-almost divisible if for each $\langle a\rangle\in \mathrm{Cu}(A)$, $k\in\mathbb{N}$ and $\varepsilon>0$, there exists $\langle b\rangle\in  \mathrm{Cu}(A)$ such that $k\langle b\rangle\leq \langle a\rangle$ and $\langle (a-\varepsilon)_+\rangle\leq(k+1)(m+1)\langle b\rangle$.
	\end{defn}
	
	\begin{defn}\cite[Definition 3.1]{RR2013}
		Let $A$ be a C*-algebra. Let $m,n\geq 1$ be integers. $A$ is said to be weakly ($m$,$n$)-divisible, if for every $\langle u\rangle\in  \mathrm{Cu}(A)$ and any $\varepsilon>0$, there exist elements $\langle x_1\rangle,\langle x_2\rangle,\cdots,\langle x_n\rangle\in  \mathrm{Cu}(A)$ such that $m\langle x_j\rangle\leq \langle u\rangle$ for all $j=1,2,\cdots,n$ and $\langle (u-\varepsilon)_+\rangle\leq \langle x_1\rangle+\cdots+\langle x_n\rangle$.
	\end{defn}
	
	Now, we will show the permanence of $m$-almost divisibility. Given an element $a\in A^\alpha$. First, we divide the proof into two cases and in each case we find a bit of extra room in the subequivalence coming from comparison. Next, we use Theorem \ref{MTF} and Theorem \ref{HMTF} to get approximately multiplicative maps $\varphi$ and $\varphi'$ and the sum of three elements to approximate $a$. Then, we can use $\varphi$ and $\varphi'$ to transfer two of them from $A$ to $A^\alpha$. The extra room is used to contral the third one. Last, we can transfer $m$-almost divisibility from $A$ to $A^\alpha$.
	\begin{thm}\label{MD}
		Let $A$ be a unital separable simple infinite dimensional stably finite C*-algebra which is $m$-almost divisible. Let $\alpha\colon G\to \mathrm{Aut}(A)$ be an action of a second-countable compact group which has the tracial Rokhlin property with comparison. Then the fixed point algebra $A^\alpha$ is $m$-almost divisible.
		\begin{proof}
			We need to show that for each $a\in (A^\alpha\otimes K)_+$, any $\varepsilon>0$, there exists $b\in (A^\alpha\otimes K)_+$ such that $k\langle b\rangle_{A^\alpha}\leq \langle a\rangle_{A^\alpha}$ and $\langle(a-\varepsilon)_+\rangle_{A^\alpha}\leq (k+1)(m+1)\langle b\rangle_{A^\alpha}$.
			
			Choose $n\in\mathbb{N}$ and $a_0\in M_n(A^\alpha)_+$ such that
			\[\|a_0-(a-\frac{\varepsilon}{3})_+\|<\frac{\varepsilon}{6}.\]
			Then we have
			\[(a-\varepsilon)_+\precsim_{A^\alpha}(a_0-\frac{\varepsilon}{3})_+\precsim_{A^\alpha}(a_0-\frac{\varepsilon}{6})_+\precsim_{A^\alpha}(a-\frac{\varepsilon}{3})_+\precsim_{A^\alpha}a.\]
			By \cite[Proposition 4.20]{MP2022}, the action $\alpha\otimes \mathrm{id}_{M_n}\colon G\to \mathrm{Aut}(A\otimes M_n)$ also has the tracial Rokhlin property with comparison. We may assume that $n=1$. 
			
			Set $a_1=(a_0-\frac{\varepsilon}{6})_+$. Since $A$ is $m$-almost divisible, and $a_1\in A^\alpha\subset A$, with any $\varepsilon'>0$, there exists $b_1\in A$ such that $k\langle b_1\rangle_A\leq \langle a_1\rangle_A$ and $\langle(a_1-\varepsilon')_+\rangle_A\leq (k+1)(m+1)\langle b_1\rangle_A$.
			
			We divide the proof into two cases.
			
			Case (1), we assume that $(a_1-\varepsilon')_+$ is Cuntz equivalent to a projection. Without loss of generality, we may assume that $(a_1-\varepsilon')_+$ is a projection.
			
			By \cite[Proposition 2.2]{PT2007}, we may assume that there exists non-zero $c\in A_+$ such that $c(a_1-\varepsilon')_+=0$ and  $\langle(a_1-\varepsilon')_+\rangle_A+ \langle c\rangle_A\leq (k+1)(m+1)\langle b_1\rangle_A$.
			
			Since $k\langle b_1\rangle_A\leq \langle a_1\rangle_A$, there exists $v=(v_{i,j})\in M_k(A)_+$ such that
			\[\|v^*diag(a_1,\underbrace{0,\cdots,0}_{k-1})v-b_1\otimes1_k\|<\frac{\varepsilon'}{2}.\]
			
			Since $\langle(a_1-\varepsilon')_+\rangle_A+\langle c\rangle_A\leq (k+1)(m+1)\langle b_1\rangle_A$, there exists $w=(w_{s,t})\in M_{(k+1)(m+1)}(A)_+$ such that
			\[\|w^*(b_1\otimes1_{(k+1)(m+1)})w-diag((a_1-\varepsilon')_++ c,\underbrace{0,\cdots,0}_{(k+1)(m+1)-1})\|<\frac{\varepsilon'}{2}.\]
			
			By Theorem \ref{MTF}, with $F_1=\{a_1,b_1,(a_1-\varepsilon')_+,c,v_{i,j},v_{i,j}^*,w_{s,t},w_{s,t}^*\colon i,j=1,2,\cdots,k\ \mathrm{and}\ s,t=1,2,\cdots,(k+1)(m+1)\}$, $F_2=\{a_1\}$, $\frac{\varepsilon'}{2((k+1)(m+1))^2}>0$ and $n=3$, there exist a projection $p\in A^\alpha$ and a unital completely positive map $\varphi\colon A\to pA^\alpha p$ such that
			\begin{enumerate}
				\item $\varphi$ is an ($3$,$F_1\cup F_2$,$\frac{\varepsilon'}{2((k+1)(m+1))^2}$)-approximately multiplicative map.
				
				\item $\|px-xp\|<\frac{\varepsilon'}{2((k+1)(m+1))^2}$ for all $x\in F_1\cup F_2$.
				
				\item $\|\varphi(x)-pxp\|<\frac{\varepsilon'}{2((k+1)(m+1))^2}$ for all $x\in F_2$.
			\end{enumerate}
			
			Thus we have
			\[\|\varphi\otimes \mathrm{id}_{M_k}(v^*)diag(\varphi(a_1),\underbrace{0,\cdots,0}_{k-1})\varphi\otimes \mathrm{id}_{M_k}(v)-\varphi(b_1)\otimes1_k\|<\varepsilon',\]
			and
			\begin{align*}
				\|\varphi\otimes &{}\mathrm{id}_{(k+1)(m+1)}(w^*)(\varphi(b_1)\otimes1_{(k+1)(m+1)})\varphi\otimes \mathrm{id}_{(k+1)(m+1)}(w)\\
				&{}-diag(\varphi((a_1-\varepsilon')_++ \varphi(c),\underbrace{0,\cdots,0}_{(k+1)(m+1)-1})\|<\varepsilon'.
			\end{align*}
			Therefore we have
			\[k\langle (\varphi(b_1)-2\varepsilon')_+\rangle_{A^\alpha}\leq \langle (\varphi(a_1)-\varepsilon')_+\rangle_{A^\alpha},\]
			and
			\[\langle (\varphi(a_1)-4\varepsilon')_+\rangle_{A^\alpha}+ \langle (\varphi(c)-3\varepsilon')_+\rangle_{A^\alpha}\leq(k+1)(m+1)\langle (\varphi(b_1)-2\varepsilon')_+\rangle_{A^\alpha}.\]
			
			Let $a_2=(1-p)a_1(1-p)$, for $a_2\in A^\alpha\subset A$, since $A$ is $m$-almost divisible, there exists $b_2\in A$ such that $k\langle b_2\rangle_A\leq \langle a_2\rangle_A$ and $\langle(a_2-\varepsilon')_+\rangle_A\leq (k+1)(m+1)\langle b_2\rangle_A$.
			
			Since $k\langle b_2\rangle_A\leq \langle a_2\rangle_A$, there exist $v'=(v'_{i,j})\in M_k(A)_+$ such that
			\[\|v'^*diag(a_2,\underbrace{0,\cdots,0}_{k-1})v'-b_2\otimes1_k\|<\varepsilon'.\]
			
			Since $\langle(a_2-\varepsilon')_+\rangle_A\leq (k+1)(m+1)\langle b_2\rangle_A$, there exists $w'=(w'_{s,t})\in M_{(k+1)(m+1)}(A)_+$ such that
			\[\|w'^*(b_2\otimes1_{(k+1)(m+1)})w'-diag((a_2-\varepsilon')_+,\underbrace{0,\cdots,0}_{(k+1)(m+1)-1})\|<\varepsilon'.\]
			
			For $A_1=(1-p)A(1-p)$, by Theorem \ref{HMTF}, with $F_1=\{b_2,a_2,(a_2-\varepsilon')_+,v_{i,j}',v_{i,j}'^*,w_{s,t}',w_{s,t}'^*\colon i,j=1,2,\cdots,k\ \mathrm{and}\ s,t=1,2,\cdots,(k+1)(m+1)\}$, $F_2=\{a_2\}$, $\frac{\varepsilon'}{2((k+1)(m+1))^2}>0$ and $n=3$, there exist a projection $q\in (A_1)^\alpha$ and a unital completely positive map $\varphi'\colon A_1\to q(A_1)^\alpha q$ such that
			\begin{itemize}
				\item[(4)] $\varphi'$ is an ($3$,$F_1\cup F_2$,$\frac{\varepsilon'}{2((k+1)(m+1))^2}$)-approximately multiplicative map.
				
				\item[(5)] $\|qx-xq\|<\frac{\varepsilon'}{2((k+1)(m+1))^2}$ for all $x\in F_1\cup F_2$.
				
				\item[(6)] $\|\varphi'(x)-qxq\|<\frac{\varepsilon'}{2((k+1)(m+1))^2}$ for all $x\in F_2$.
				
				\item[(7)] $1_{A_1}-q\precsim_{(A_1)^\alpha} (\varphi(c)-3\varepsilon')_+$. 
			\end{itemize}
			
			Thus we have
			\[\|\varphi'\otimes \mathrm{id}_{M_k}(v'^*)diag(\varphi'(a_2),\underbrace{0,\cdots,0}_{k-1})\varphi'\otimes \mathrm{id}_{M_k}(v')-\varphi'(b_2)\otimes1_k\|<\varepsilon',\]
			and
			\begin{align*}
				\|\varphi'\otimes &{}\mathrm{id}_{(k+1)(m+1)}(w'^*)(\varphi'(b_2)\otimes1_{(k+1)(m+1)})\varphi'\otimes \mathrm{id}_{(k+1)(m+1)}(w')\\
				&{}-diag(\varphi'((a_2-\varepsilon')_+,\underbrace{0,\cdots,0}_{(k+1)(m+1)-1})\|<\varepsilon'.
			\end{align*}
			Therefore we have
			\[k\langle (\varphi'(b_2)-4\varepsilon')_+\rangle_{A^\alpha}\leq \langle (\varphi'(a_2)-2\varepsilon')_+\rangle_{A^\alpha},\]
			and
			\[\langle (\varphi'(a_2)-6\varepsilon')_+\rangle_{A^\alpha}\leq(k+1)(m+1)\langle (\varphi'(b_2)-4\varepsilon')_+\rangle_{A^\alpha}.\]
			
			Therefore, with $\varepsilon'$ sufficiently small, we have
			\begin{align*}
				&{}k\langle(\varphi(b_1)-2\varepsilon')_+\oplus(\varphi'(b_2)-4\varepsilon')_+\rangle_{A^\alpha}\\
				=&{}k\langle(\varphi(b_1)-2\varepsilon')_+\rangle_{A^\alpha}+ k\langle(\varphi'(b_2)-4\varepsilon')_+\rangle_{A^\alpha}\\
				\leq&{}\langle(\varphi(a_1)-\varepsilon')_+\rangle_{A^\alpha}+\langle(\varphi'(a_2)-2\varepsilon')_+\rangle_{A^\alpha}\\
				\leq&{}\langle(\varphi(a_1)-\varepsilon')_+\rangle_{A^\alpha}+\langle(\varphi'(a_2)-2\varepsilon')_+\rangle_{A^\alpha}+\langle((1_{A_1}-q)a_2(1_{A_1}-q)-2\varepsilon')_+\rangle_{A^\alpha}\\
				=&{}\langle (\varphi(a_1)-\varepsilon')_+\oplus (\varphi'(a_2)-2\varepsilon')_+\oplus ((1_{A_1}-q)a_2(1_{A_1}-q)-2\varepsilon')_+ \rangle_{A^\alpha}\\
				=&{}\langle (\varphi(a_1)-\varepsilon')_++ (\varphi'(a_2)-2\varepsilon')_++ ((1_{A_1}-q)a_2(1_{A_1}-q)-2\varepsilon')_+ \rangle_{A^\alpha}\\
				\leq&{}\langle (\varphi(a_1)-\varepsilon')_++ ((1-p)a_1(1-p)-\varepsilon')_+\leq\langle a_1\rangle_{A^\alpha}=\langle (a_0-\frac{\varepsilon}{6})_+\rangle_{A^\alpha},
			\end{align*}
			and we also have
			\begin{align*}
				&{}\langle (a_0-\frac{\varepsilon}{3})_+\rangle_{A^\alpha}=\langle (a_1-\frac{\varepsilon}{6})_+\rangle_{A^\alpha}\\
				\leq&{}\langle (\varphi(a_1)-4\varepsilon')_++(\varphi'(a_2)-6\varepsilon')_++((1_{A_1}-q)a_2(1_{A_1}-q)-6\varepsilon')_+\rangle_{A^\alpha}\\
				=&{}\langle (\varphi(a_1)-4\varepsilon')_+\oplus(\varphi'(a_2)-6\varepsilon')_+\oplus((1_{A_1}-q)a_2(1_{A_1}-q)-6\varepsilon')_+\rangle_{A^\alpha}\\
				=&{}\langle (\varphi(a_1)-4\varepsilon')_+\rangle_{A^\alpha}+\langle(\varphi'(a_2)-6\varepsilon')_+\rangle_{A^\alpha}+\langle((1_{A_1}-q)a_2(1_{A_1}-q)-6\varepsilon')_+\rangle_{A^\alpha}\\
				\leq&{}\langle(\varphi(a_1)-4\varepsilon')_+\rangle_{A^\alpha}+\langle(\varphi'(a_2)-6\varepsilon')_+\rangle_{A^\alpha}+\langle1_{A_1}-q\rangle_{A^\alpha}\\
				\leq&{}\langle(\varphi(a_1)-4\varepsilon')_+\rangle_{A^\alpha}+\langle(\varphi'(a_2)-6\varepsilon')_+\rangle_{A^\alpha}+\langle (\varphi(c)-3\varepsilon')_+\rangle_{A^\alpha}\\
				\leq&{}(k+1)(m+1)\langle(\varphi(b_1)-2\varepsilon')_+\rangle_{A^\alpha}+(k+1)(m+1)\langle(\varphi'(b_2)-4\varepsilon')_+\rangle_{A^\alpha}\\
				=&{}(k+1)(m+1)\langle(\varphi(b_1)-2\varepsilon')_+\oplus(\varphi'(b_2)-4\varepsilon')_+\rangle_{A^\alpha}.
			\end{align*}
			
			Set $b=(\varphi(b_1)-2\varepsilon')_+\oplus(\varphi'(b_2)-4\varepsilon')_+$, we have
			\[k\langle b\rangle_{A^\alpha}\leq \langle (a_0-\frac{\varepsilon}{6})_+\rangle_{A^\alpha}\leq \langle a\rangle_{A^\alpha}\ \mathrm{and}\ \langle (a-\varepsilon)_+\rangle_{A^\alpha}\leq(a_0-\frac{\varepsilon}{3})_+\rangle_{A^\alpha}\leq(k+1)(m+1)\langle b\rangle_{A^\alpha}.\]
			
			Case (2), we assume that $(a_1-\varepsilon')_+$ is not Cuntz equivalent to a projection.
			
			By \cite[Theorem 2.1(4)]{EFF2018}, there is a non-zero positive element $d\in A_+$ such that $d(a_1-2\varepsilon')_+=0$ and $\langle(a_1-2\varepsilon')_+\rangle_A+\langle d\rangle_A\leq\langle(a_1-\varepsilon')_+\rangle_A$.
			
			Since $k\langle b_1\rangle_A\leq \langle a_1\rangle_A$, there exists $v=(v_{i,j})\in M_k(A)_+$ such that
			\[\|v^*diag(a,\underbrace{0,\cdots,0}_{k-1})v-b_1\otimes1_k\|<\frac{\varepsilon'}{2}.\]
			
			Since $\langle(a_1-2\varepsilon')_+\rangle_A+\langle d\rangle_A\leq\langle(a_1-\varepsilon')_+\rangle_A\leq(k+1)(m+1)\langle b_1\rangle_A$, there exists $w=(w_{s,t})\in M_{(k+1)(m+1)}(A)_+$ such that
			\[\|w^*(b_1\otimes1_{(k+1)(m+1)})w-diag((a_1-2\varepsilon')_++ d,\underbrace{0,\cdots,0}_{(k+1)(m+1)-1})\|<\frac{\varepsilon'}{2}.\]
			
			By Theorem \ref{MTF}, with $F_1=\{a_1,b_1,(a_1-2\varepsilon')_+,d,v_{i,j},v_{i,j}^*,w_{s,t},w_{s,t}^*\colon i,j=1,2,\cdots,k\ \mathrm{and}\ s,t=1,2,\cdots,(k+1)(m+1)\}$, $F_2=\{a\}$, $\frac{\varepsilon'}{2((k+1)(m+1))^2}>0$ and $n=3$, there exist a projection $p\in A^\alpha$ and a unital completely positive map $\varphi\colon A\to pA^\alpha p$ such that
			\begin{itemize}
				\item[(8)] $\varphi$ is an ($3$,$F_1\cup F_2$,$\frac{\varepsilon'}{2((k+1)(m+1))^2}$)-approximately multiplicative map.
				
				\item[(9)] $\|px-xp\|<\frac{\varepsilon'}{2((k+1)(m+1))^2}$ for all $x\in F_1\cup F_2$.
				
				\item[(10)] $\|\varphi(x)-pxp\|<\frac{\varepsilon'}{2((k+1)(m+1))^2}$ for all $x\in F_2$.
			\end{itemize}
			
			Thus we have
			\[\|\varphi\otimes \mathrm{id}_{M_k}(v^*)diag(\varphi(a_1),\underbrace{0,\cdots,0}_{k-1})\varphi\otimes \mathrm{id}_{M_k}(v)-\varphi(b_1)\otimes1_k\|<\varepsilon',\]
			and
			\begin{align*}
				\|\varphi\otimes &{}\mathrm{id}_{(k+1)(m+1)}(w^*)(\varphi(b_1)\otimes1_{(k+1)(m+1)})\varphi\otimes \mathrm{id}_{(k+1)(m+1)}(w)\\
				&{}-diag(\varphi((a_1-2\varepsilon')_++ d,\underbrace{0,\cdots,0}_{(k+1)(m+1)-1})\|<\varepsilon'.
			\end{align*}
			Therefore we have
			\[k\langle (\varphi(b_1)-2\varepsilon')_+\rangle_{A^\alpha}\leq \langle(\varphi(a_1)-\varepsilon')_+\rangle_{A^\alpha},\]
			and
			\[\langle(\varphi(a_1)-5\varepsilon')_+\rangle_{A^\alpha}+\langle(\varphi(d)-3\varepsilon)_+\rangle_{A^\alpha}\leq(k+1)(m+1)\langle(\varphi(b_1)-2\varepsilon')_+\rangle_{A^\alpha}.\]
			
			Let $a_2=(1-p)a_1(1-p)$, for $a_2\in A^\alpha\subset A$, since $A$ is $m$-almost divisible, there exists $b_2\in A$ such that $k\langle b_2\rangle_A\leq a_2$ and $\langle(a_2-\varepsilon')_+\rangle_A\leq (k+1)(m+1)b_2$.
			
			Since $k\langle b_2\rangle_A\leq \langle a_2\rangle_A$, there exist $v'=(v'_{i,j})\in M_k(A)_+$ such that
			\[\|v'^*diag(a_2,\underbrace{0,\cdots,0}_{k-1})v'-b_2\otimes1_k\|<\frac{\varepsilon'}{2},\]
			
			Since $\langle(a_2-\varepsilon')_+\rangle_A\leq (k+1)(m+1)\langle b_2\rangle_A$, there exists $w'=(w'_{s,t})\in M_{(k+1)(m+1)}(A)_+$ such that
			\[\|w'^*(b_2\otimes1_{(k+1)(m+1)})w'-diag((a_2-\varepsilon')_+,\underbrace{0,\cdots,0}_{(k+1)(m+1)-1})\|<\frac{\varepsilon'}{2}.\]
			
			For $A_1=(1-p)A(1-p)$, by Theorem \ref{HMTF}, with $F_1=\{b_2,a_2,(a_1-\varepsilon')_+,v_{i,j}',v_{i,j}'^*,w_{s,t}',w_{s,t}'^*\colon i,j=1,2,\cdots,k\ \mathrm{and}\ s,t=1,2,\cdots,(k+1)(m+1)\}$, $F_2=\{a_2\}$, $\frac{\varepsilon'}{2((k+1)(m+1))^2}>0$ and $n=3$, there exist a projection $q\in (A_1)^\alpha$ and a unital completely positive map $\varphi'\colon A_1\to q(A_1)^\alpha q$ such that
			\begin{itemize}
				\item[(11)] $\varphi'$ is an ($3$,$F_1\cup F_2$,$\frac{\varepsilon'}{2((k+1)(m+1))^2}$)-approximately multiplicative map.
				
				\item[(12)] $\|qx-xq\|<\frac{\varepsilon'}{2((k+1)(m+1))^2}$ for all $x\in F_1\cup F_2$.
				
				\item[(13)] $\|\varphi'(x)-qxq\|<\frac{\varepsilon'}{2((k+1)(m+1))^2}$ for all $x\in F_2$.
				
				\item[(14)] $1_{A_1}-q\precsim_{(A_1)^\alpha} (\varphi(d)-3\varepsilon')_+$. 
			\end{itemize}
			
			Thus we have
			\[\|\varphi'\otimes \mathrm{id}_{M_k}(v'^*)diag(\varphi'(a_2),\underbrace{0,\cdots,0}_{k-1})\varphi'\otimes \mathrm{id}_{M_k}(v')-\varphi'(b_2)\otimes1_k\|<\varepsilon',\]
			and
			\begin{align*}
				\|\varphi'\otimes &{}\mathrm{id}_{(k+1)(m+1)}(w'^*)(\varphi'(b_2)\otimes1_{(k+1)(m+1)})\varphi'\otimes \mathrm{id}_{(k+1)(m+1)}(w')\\
				&{}-diag(\varphi'((a_2-\varepsilon')_+,\underbrace{0,\cdots,0}_{(k+1)(m+1)-1})\|<\varepsilon'.
			\end{align*}
			Therefore we have
			\[k\langle(\varphi'(b_2)-4\varepsilon')_+\rangle_{A^\alpha}\leq \langle(\varphi'(a_2)-2\varepsilon')_+\rangle_{A^\alpha},\]
			and
			\[\langle(\varphi'(a_2)-6\varepsilon')_+\rangle_{A^\alpha}\leq(k+1)(m+1)\langle(\varphi'(b_2)-4\varepsilon')_+\rangle_{A^\alpha}.\]
			
			Therefore, with $\varepsilon'$ sufficiently small, we have
			\begin{align*}
				&{}k\langle(\varphi(b_1)-2\varepsilon')_+\oplus(\varphi'(b_2)-4\varepsilon')_+\rangle_{A^\alpha}\\
				=&{}k\langle(\varphi(b_1)-2\varepsilon')_+\rangle_{A^\alpha}+ k\langle(\varphi'(b_2)-4\varepsilon')_+\rangle_{A^\alpha}\\
				\leq&{}\langle(\varphi(a_1)-\varepsilon')_+\rangle_{A^\alpha}+\langle(\varphi'(a_2)-2\varepsilon')_+\rangle_{A^\alpha}\\
				\leq&{}\langle(\varphi(a_1)-\varepsilon')_+\rangle_{A^\alpha}+\langle(\varphi'(a_2)-2\varepsilon')_+\rangle_{A^\alpha}+\langle((1_{A_1}-q)a_2(1_{A_1}-q)-2\varepsilon')_+\rangle_{A^\alpha}\\
				=&{}\langle (\varphi(a_1)-\varepsilon')_+\oplus (\varphi'(a_2)-2\varepsilon')_+\oplus ((1_{A_1}-q)a_2(1_{A_1}-q)-2\varepsilon')_+ \rangle_{A^\alpha}\\
				=&{}\langle (\varphi(a_1)-\varepsilon')_++ (\varphi'(a_2)-2\varepsilon')_++ ((1_{A_1}-q)a_2(1_{A_1}-q)-2\varepsilon')_+ \rangle_{A^\alpha}\\
				\leq&{}\langle (\varphi(a_1)-\varepsilon')_++ ((1-p)a_1(1-p)-\varepsilon')_+\leq\langle a_1\rangle_{A^\alpha}=\langle (a_0-\frac{\varepsilon}{6})_+\rangle_{A^\alpha},
			\end{align*}
			and we also have
			\begin{align*}
				&{}\langle (a_0-\frac{\varepsilon}{3})_+\rangle_{A^\alpha}=\langle (a_1-\frac{\varepsilon}{6})_+\rangle_{A^\alpha}\\
				\leq&{}\langle (\varphi(a_1)-5\varepsilon')_++(\varphi'(a_2)-6\varepsilon')_++((1_{A_1}-q)a_2(1_{A_1}-q)-6\varepsilon')_+\rangle_{A^\alpha}\\
				=&{}\langle (\varphi(a_1)-5\varepsilon')_+\oplus(\varphi'(a_2)-6\varepsilon')_+\oplus((1_{A_1}-q)a_2(1_{A_1}-q)-6\varepsilon')_+\rangle_{A^\alpha}\\
				=&{}\langle (\varphi(a_1)-5\varepsilon')_+\rangle_{A^\alpha}+\langle(\varphi'(a_2)-6\varepsilon')_+\rangle_{A^\alpha}+\langle((1_{A_1}-q)a_2(1_{A_1}-q)-6\varepsilon')_+\rangle_{A^\alpha}\\
				\leq&{}\langle(\varphi(a_1)-5\varepsilon')_+\rangle_{A^\alpha}+\langle(\varphi'(a_2)-6\varepsilon')_+\rangle_{A^\alpha}+\langle1_{A_1}-q\rangle_{A^\alpha}\\
				\leq&{}\langle(\varphi(a_1)-5\varepsilon')_+\rangle_{A^\alpha}+\langle(\varphi'(a_2)-6\varepsilon')_+\rangle_{A^\alpha}+\langle (\varphi(d)-3\varepsilon')_+\rangle_{A^\alpha}\\
				\leq&{}(k+1)(m+1)\langle(\varphi(b_1)-2\varepsilon')_+\rangle_{A^\alpha}+(k+1)(m+1)\langle(\varphi'(b_2)-4\varepsilon')_+\rangle_{A^\alpha}\\
				=&{}(k+1)(m+1)\langle(\varphi(b_1)-2\varepsilon')_+\oplus(\varphi'(b_2)-4\varepsilon')_+\rangle_{A^\alpha}.
			\end{align*}
			
			Set $b=(\varphi(b_1)-2\varepsilon')_+\oplus(\varphi'(b_2)-4\varepsilon')_+$, we have
			\[k\langle b\rangle_{A^\alpha}\leq \langle (a_0-\frac{\varepsilon}{6})_+\rangle_{A^\alpha}\leq \langle a\rangle_{A^\alpha}\ \mathrm{and}\ \langle (a-\varepsilon)_+\rangle_{A^\alpha}\leq(a_0-\frac{\varepsilon}{3})_+\rangle_{A^\alpha}\leq(k+1)(m+1)\langle b\rangle_{A^\alpha}.\]
		\end{proof}
	\end{thm}
	
	\begin{cor}
		Let $A$ be a unital separable simple infinite dimensional stably finite C*-algebra which is $m$-almost divisible. Let $\alpha\colon G\to \mathrm{Aut}(A)$ be an action of a second-countable compact group which has the tracial Rokhlin property with comparison. Then the crossed product $A\rtimes_\alpha G$ is $m$-almost divisible.
		\begin{proof}
			The algebra $A^\alpha$ is $m$-almost divisible by Theorem \ref{MD}. Theorem \ref{CPS} implies that $A\rtimes_\alpha G$ is stably isomorphic to $A^\alpha$. Since two stably isomorphic C*-algebras have canonically isomorphic Cuntz semigroups, we know that $A\rtimes_\alpha G$ is $m$-almost divisible.
		\end{proof}
	\end{cor}
	
	Now, we will show the permanence of  weakly ($m$,$n$)-divisibility. The idea of the proof is essentially the same as that of Theorem \ref{MD}. 
	\begin{thm}\label{WMND}
		Let $A$ be a unital separable simple infinite dimensional stably finite C*-algebra which is weakly ($m$,$n$)-divisible. Let $\alpha\colon G\to \mathrm{Aut}(A)$ be an action of a second-countable compact group which has the tracial Rokhlin property with comparison. Then the fixed point algebra $A^\alpha$ is weakly ($m$,$n$)-divisible.
		\begin{proof}
			We need to show that for each $a\in (A^\alpha\otimes K)_+$, any $\varepsilon>0$, there exists $x_1,x_2,\cdots,x_n\in (A^\alpha\otimes K)_+$ such that $m\langle x_j\rangle_{A^\alpha}\leq \langle a\rangle_{A^\alpha}$ for all $1\leq j\leq n$, and $\langle(a-\varepsilon)_+\rangle_{A^\alpha}\leq \langle\oplus_{i=1}^n x_i\rangle_{A^\alpha}$. 
			
			Choose $n\in\mathbb{N}$ and $a_0\in M_n(A^\alpha)_+$ such that
			\[\|a_0-(a-\frac{\varepsilon}{3})_+\|<\frac{\varepsilon}{6}.\]
			Then we have
			\[(a-\varepsilon)_+\precsim_{A^\alpha}(a_0-\frac{\varepsilon}{3})_+\precsim_{A^\alpha}(a_0-\frac{\varepsilon}{6})_+\precsim_{A^\alpha}(a-\frac{\varepsilon}{3})_+\precsim_{A^\alpha}a.\]
			By \cite[Proposition 4.20]{MP2022}, the action $\alpha\otimes \mathrm{id}_{M_n}\colon G\to \mathrm{Aut}(A\otimes M_n)$ also has the tracial Rokhlin property with comparison. We may assume that $n=1$. 
			
			Set $a_1=(a_0-\frac{\varepsilon}{6})_+$. Since $A$ is weakly ($m$,$n$)-divisible, and $a_1\in A^\alpha\subset A$, with any $\varepsilon'>0$, there exists $y_1,y_2,\cdots,y_n\in A$ such that $m\langle y_j\rangle_A\leq \langle a_1\rangle_A$ for all $1\leq j\leq n$, and $\langle(a_1-\varepsilon')_+\rangle_A\leq \langle\oplus_{i=1}^n y_i\rangle_A$.
			
			We divide the proof into two case.
			
			Case (1), we assume that $(a_1-\varepsilon')_+$ is Cuntz equivalent to a projection. Without loss of generality, we may assume that $(a_1-\varepsilon')_+$ is a projection.
			
			As in the same step in the arguement in Theorem \ref{MD}, we get a projection $p\in A^\alpha$ and a unital completely positive map $\varphi\colon A\to pA^\alpha p$ such that
			\[m\langle (\varphi(y_j)-2\varepsilon')_+\rangle_{A^\alpha}\leq \langle (\varphi(a_1)-\varepsilon')_+\rangle_{A^\alpha},\]
			for all $1\leq j\leq n$, and
			\begin{align*}
				\langle (\varphi(a_1)-4\varepsilon')_+\rangle_{A^\alpha}+ \langle (\varphi(c)-3\varepsilon')_+\rangle_{A^\alpha}
				\leq\langle (\varphi(\oplus_{i=1}^n y_i)-2\varepsilon')_+\rangle_{A^\alpha}
				=\langle \oplus_{i=1}^n (\varphi(y_i)-2\varepsilon')_+\rangle_{A^\alpha}.
			\end{align*}
			
			Let $a_2=(1-p)a_1(1-p)$, for $a_2\in A^\alpha\subset A$, since $A$ is weakly ($m$,$n$)-divisible, there exists $z_1,z_2,\cdots,z_n\in A$ such that $m\langle z_j\rangle_A\leq \langle a_2\rangle_A$ for all $1\leq j\leq n$ and $\langle(a_2-\varepsilon')_+\rangle_A\leq \langle \oplus_{i=1}^n z_i\rangle_A$.
			
			As in the same step in the arguement in Theorem \ref{MD}, we get a projection $p\in A^\alpha$ and a unital completely positive map $\varphi'\colon A_1\to q(A_1)^\alpha q$ such that $1_{A_1}-q\precsim_{(A_1)^\alpha} (\varphi(c)-3\varepsilon')_+$ and 
			\[m\langle (\varphi'(z_j)-4\varepsilon')_+\rangle_{A^\alpha}\leq \langle (\varphi'(a_2)-2\varepsilon')_+\rangle_{A^\alpha}.\]
			for all $1\leq j\leq n$, and
			\begin{align*}
				\langle (\varphi'(a_2)-6\varepsilon')_+\rangle_{A^\alpha}
				\leq\langle (\varphi'(\oplus_{i=1}^n z_i)-4\varepsilon')_+\rangle_{A^\alpha}
				=\langle \oplus_{i=1}^n (\varphi'(z_i)-4\varepsilon')_+\rangle_{A^\alpha}.
			\end{align*}
			
			Therefore, with $\varepsilon'$ sufficiently small, as in the same step in the arguement in Theorem \ref{MD}, we have
			\begin{align*}
				m\langle(\varphi(y_j)-2\varepsilon')_+\oplus(\varphi'(z_j)-4\varepsilon')_+\rangle_{A^\alpha}
				\leq\langle a_1\rangle_{A^\alpha}=\langle (a_0-\frac{\varepsilon}{6})_+\rangle_{A^\alpha},
			\end{align*}
			for all $1\leq j\leq n$, and
			\begin{align*}
				\langle (a_0-\frac{\varepsilon}{3})_+\rangle_{A^\alpha}=\langle (a_1-\frac{\varepsilon}{6})_+\rangle_{A^\alpha}
				\leq\langle \oplus_{i=1}^n ((\varphi(y_i)-2\varepsilon')_+\oplus(\varphi'(z_i)-4\varepsilon')_+)\rangle_{A^\alpha}.
			\end{align*}
			
			Set $x_j=(\varphi(y_j)-2\varepsilon')_+\oplus(\varphi'(z_j)-4\varepsilon')_+$ for all $1\leq j\leq n$, we have
			\[m\langle x_j\rangle_{A^\alpha}\leq \langle (a_0-\frac{\varepsilon}{6})_+\rangle_{A^\alpha}\leq \langle a\rangle_{A^\alpha}\ \mathrm{and}\ \langle (a-\varepsilon)_+\rangle_{A^\alpha}\leq(a_0-\frac{\varepsilon}{3})_+\rangle_{A^\alpha}\leq\langle \oplus_{i=1}^n x_i\rangle_{A^\alpha}.\]
			
			Case (2), we assume that $(a-\varepsilon')_+$ is not Cuntz equivalent to a projection.
			
			By \cite[Theorem 2.1(4)]{EFF2018}, there is a non-zero positive element $d\in A_+$ such that $d(a-2\varepsilon')_+=0$ and $\langle(a-2\varepsilon')_+\rangle_A+\langle d\rangle_A\leq\langle(a-\varepsilon')_+\rangle_A$.
			
			As in the same step in the arguement in Theorem \ref{MD}, we get a projection $p\in A^\alpha$ and a unital completely positive map $\varphi\colon A\to pA^\alpha p$ such that
			\[m\langle (\varphi(y_j)-2\varepsilon')_+\rangle_{A^\alpha}\leq \langle (\varphi(a)-\varepsilon')_+\rangle_{A^\alpha},\]
			for all $1\leq j\leq n$, and
			\begin{align*}
				\langle (\varphi(a)-5\varepsilon')_+\rangle_{A^\alpha}+ \langle (\varphi(d)-3\varepsilon')_+\rangle_{A^\alpha}
				\leq\langle (\varphi(\oplus_{i=1}^n y_i)-2\varepsilon')_+\rangle_{A^\alpha}
				=\langle \oplus_{i=1}^n (\varphi(y_i)-2\varepsilon')_+\rangle_{A^\alpha}.
			\end{align*}
			
			As in the same step in the arguement in Theorem \ref{MD}, we get a projection $p\in A^\alpha$ and a unital completely positive map $\varphi'\colon A_1\to q(A_1)^\alpha q$ such that $1_{A_1}-q\precsim_{(A_1)^\alpha} (\varphi(d)-3\varepsilon')_+$ and
			\[m\langle (\varphi'(z_j)-4\varepsilon')_+\rangle_{A^\alpha}\leq \langle (\varphi'(a_1)-2\varepsilon')_+\rangle_{A^\alpha}.\]
			for all $1\leq j\leq n$, and
			\begin{align*}
				\langle (\varphi'(a_1)-6\varepsilon')_+\rangle_{A^\alpha}
				\leq\langle (\varphi'(\oplus_{i=1}^n z_i)-4\varepsilon')_+\rangle_{A^\alpha}
				=\langle \oplus_{i=1}^n (\varphi'(z_i)-4\varepsilon')_+\rangle_{A^\alpha}.
			\end{align*}
			
			Therefore, with $\varepsilon'$ sufficiently small, as in the same step in the arguement in Theorem \ref{MD}, we have
			\begin{align*}
				m\langle(\varphi(y_j)-2\varepsilon')_+\oplus(\varphi'(z_j)-4\varepsilon')_+\rangle_{A^\alpha}
				\leq\langle a_1\rangle_{A^\alpha}=\langle (a_0-\frac{\varepsilon}{6})_+\rangle_{A^\alpha},
			\end{align*}
			for all $1\leq j\leq n$, and
			\begin{align*}
				\langle (a_0-\frac{\varepsilon}{3})_+\rangle_{A^\alpha}=\langle (a_1-\frac{\varepsilon}{6})_+\rangle_{A^\alpha}
				\leq\langle \oplus_{i=1}^n ((\varphi(y_i)-2\varepsilon')_+\oplus(\varphi'(z_i)-4\varepsilon')_+)\rangle_{A^\alpha}.
			\end{align*}
			
			Set $x_j=(\varphi(y_j)-2\varepsilon')_+\oplus(\varphi'(z_j)-4\varepsilon')_+$ for all $1\leq j\leq n$, we have
			\[m\langle x_j\rangle_{A^\alpha}\leq \langle (a_0-\frac{\varepsilon}{6})_+\rangle_{A^\alpha}\leq \langle a\rangle_{A^\alpha}\ \mathrm{and}\ \langle (a-\varepsilon)_+\rangle_{A^\alpha}\leq(a_0-\frac{\varepsilon}{3})_+\rangle_{A^\alpha}\leq\langle \oplus_{i=1}^n x_i\rangle_{A^\alpha}.\]
		\end{proof}
	\end{thm}
	
	\begin{cor}
		Let $A$ be a unital separable simple infinite dimensional stably finite C*-algebra which is weakly ($m$,$n$)-divisible. Let $\alpha:G\rightarrow Aut(A)$ be an action of a second-countable compact group which has the tracial Rokhlin property with comparison. Then the crossed product $A\rtimes_\alpha G$ is weakly ($m$,$n$)-divisible.
		\begin{proof}
			The algebra $A^\alpha$ is weakly ($m$,$n$)-divisible by Theorem \ref{WMND}. Theorem \ref{CPS} implies that $A\rtimes_\alpha G$ is  stably isomorphic to $A^\alpha$. Since two stably isomorphic $C^*$-algebras have canonically isomorphic Cuntz semigroups, we know that $A\rtimes_\alpha G$ is weakly ($m$,$n$)-divisible.
		\end{proof}
	\end{cor}
	
	\section{Futher work}
	Gardella and Santiago introduced the equivariant Cuntz semigroups for compact group actions on C*-algebras in \cite{GS2017}. They proved that the equivariant Cuntz semigroup, as a functor, is continuous and stable. One of their main results is an analog of Julg's theorem as following.
	\begin{thm}\cite[Theorem 5.14]{GS2017}
		Let $G$ be a compact group, let $A$ be a C*-algebra, and let $\alpha\colon G\to \mathrm{Aut}(A)$ be a continuous action. Then there exists a natural $\mathrm{Cu}^G$-isomorphism
		\[\mathrm{Cu}^G(A,\alpha)\cong \mathrm{Cu}(A\rtimes_\alpha G),\]
		where the $\mathrm{Cu}(G)$-semimodule structure on $\mathrm{Cu}(A\rtimes_\alpha G)$ is given by \cite[Definition 5.10]{GS2017}.
	\end{thm}
	Moreover, they computed the equivariant Cuntz semigroups of a number of dynamical systems, particularly, the Rokhlin property for an action of a finite group and got the following result.
	\begin{prop}\cite[Proposition 6.2]{GS2017}
		Let $G$ be a finite group, let $A$ be a C*-algebra, and let $\alpha\colon G\to \mathrm{Aut}(A)$ be an action with the Rokhlin property. Then there exists a natural $\mathrm{Cu}^G$-isomorphism
		\[\mathrm{Cu}^G(A,\alpha)\cong \mathrm{Cu}(A)^{\mathrm{Cu}(\alpha)}_\mathbb{N},\]
		where the the induced $\mathrm{Cu}(G)$-semimodule structure on $\mathrm{Cu}(A)^{\mathrm{Cu}(\alpha)}_\mathbb{N}$ is trivial.
	\end{prop}
	For the compact group action with the Rokhlin property, this is relatively easy, and we imagine that the outcome for the tracial Rokhlin property (with comparion) will be much more interesting since the $\mathrm{Cu}(G)$-semimodule structure will not be trivial in general. We will try to compute it in the future work.
	
	\subsection*{Acknowledgments}
	The authors would like to thank the referees for their helpful comments. This work is partially supported by National Natural Science Foundation of China [Grant no. 11871375].

\end{document}